\documentclass[11pt]{article}
\usepackage{amsmath}
\usepackage{amssymb}
\usepackage{amsthm}
\usepackage{latexsym}
\usepackage[dvips]{graphicx,color}
\usepackage{float}
\usepackage{array}
\usepackage{xcolor}
\usepackage{enumerate} 
\usepackage{titlesec}
\usepackage{multirow}
\usepackage[hyphens]{url}
\usepackage{hyperref}
\usepackage{etoolbox}
\usepackage{MnSymbol}
\usepackage{ifthen}
\usepackage{mathrsfs}
\usepackage{xspace}
\usepackage[normalem]{ulem}

\newtheoremstyle{mystyle}%
{3pt}
{3pt}
{\color{blue}}
{}
{\bfseries\color{blue}}
{.}
{.5em}
{}

\theoremstyle{plain}
\newtheorem{theorem}{Theorem}[section]

\newtheorem{proposition}[theorem]{Proposition}
\newtheorem{corollary}[theorem]{Corollary}

\newtheorem{lemma}[theorem]{Lemma}

\theoremstyle{definition}

\theoremstyle{remark}
\newtheorem{remark}[theorem]{Remark}

\theoremstyle{mystyle}

\newcommand{\fequiv}[1]{\ensuremath{\equiv_{#1, \fo}}}
\newcommand{\mequiv}[1]{\ensuremath{\equiv_{#1, \mso}}}
\newcommand{\lequiv}[1]{\ensuremath{\equiv_{#1, \mc{L}}}}

\newcommand{\lth}[2]{\ensuremath{\mathsf{Th}_{#1,
      \mc{L}}(#2)}}
\newcommand{\msoth}[2]{\ensuremath{\mathsf{Th}_{#1,
      \mso}(#2)}}
\newcommand{\foth}{\ensuremath{\mathsf{Th}}}
\renewcommand{\foth}[1]{\ensuremath{\mathrm{FO}\text{-}\mathsf{Th}(#1)}}
\renewcommand{\msoth}[1]{\ensuremath{\mathrm{MSO}\text{-}\mathsf{Th}(#1)}}
\renewcommand{\lth}[1]{\ensuremath{\mc{L}\text{-}\mathsf{Th}(#1)}}

\newcommand{\fo}{\ensuremath{\text{FO}}}
\newcommand{\mso}{\ensuremath{\text{MSO}}}

\newcommand{\rank}[1]{\ensuremath{\text{rank}(#1)}}

\newcommand{\lt}{{{\L}o{\'s}-Tarski}}

\newcommand{\eat}[1]{}

\newcommand{\tree}[1]{\ensuremath{\mathsf{#1}}}
\newcommand{\forest}[1]{\ensuremath{\mathsf{#1}}}

\newcommand{\cl}[1]{\ensuremath{\mathcal{#1}}}

\newcommand{\str}[1]{\ensuremath{\mathsf{Str}(#1)}}

\newcommand{\tbf}[1]{\textbf{#1}}
\newcommand{\mc}[1]{\mathcal{#1}}

\newcommand{\mf}[1]{\mathfrak{#1}}

\makeatletter
\newsavebox\myboxA
\newsavebox\myboxB
\newlength\mylenA

\newcommand*\xoverline[2][0.75]{%
    \sbox{\myboxA}{$\m@th#2$}%
    \setbox\myboxB\null
    \ht\myboxB=\ht\myboxA%
    \dp\myboxB=\dp\myboxA%
    \wd\myboxB=#1\wd\myboxA
    \sbox\myboxB{$\m@th\overline{\copy\myboxB}$}
    \setlength\mylenA{\the\wd\myboxA}
    \addtolength\mylenA{-\the\wd\myboxB}%
    \ifdim\wd\myboxB<\wd\myboxA%
       \rlap{\hskip 0.5\mylenA\usebox\myboxB}{\usebox\myboxA}%
    \else
        \hskip -0.5\mylenA\rlap{\usebox\myboxA}{\hskip 0.5\mylenA\usebox\myboxB}%
    \fi}
    
  \newcommand*\xunderline[2][0.75]{%
    \sbox{\myboxA}{$\m@th#2$}%
    \setbox\myboxB\null
    \ht\myboxB=\ht\myboxA%
    \dp\myboxB=\dp\myboxA%
    \wd\myboxB=#1\wd\myboxA
    \sbox\myboxB{$\m@th\underline{\copy\myboxB}$}
    \setlength\mylenA{\the\wd\myboxA}
    \addtolength\mylenA{-\the\wd\myboxB}%
    \ifdim\wd\myboxB<\wd\myboxA%
       \rlap{\hskip 0.5\mylenA\usebox\myboxB}{\usebox\myboxA}%
    \else
        \hskip -0.5\mylenA\rlap{\usebox\myboxA}{\hskip 0.5\mylenA\usebox\myboxB}%
    \fi}
\makeatother

\renewcommand{\str}[1]{\ensuremath{\mathcal{#1}}}
\newcommand{\tower}[1]{\ensuremath{\mathsf{tower}(#1)}}

\renewcommand{\forest}[1]{\ensuremath{\mathsf{#1}}}
\renewcommand{\root}[1]{\ensuremath{\mathsf{root}(#1)}}
\newcommand{\rootref}{\ensuremath{\mathsf{root}}}
\newcommand{\mdelta}[2]{\ensuremath{\delta_{#2}(#1)}}

\renewcommand{\mdelta}[2]{\ensuremath{\Delta^{\mso}_{#1, #2}}}

\newcommand{\ef}{Ehrenfeucht-Fr\"aiss\'e}

\renewcommand{\mso}{\ensuremath{\mathrm{MSO}}}
\renewcommand{\fo}{\ensuremath{\mathrm{FO}}}

\renewcommand{\rank}[1]{\ensuremath{\mbox{rank}(#1)}}
\renewcommand{\lequiv}[1]{\ensuremath{\equiv_{#1, \mc{L}}}}

\newcommand{\ind}[2]{\ensuremath{\mathcal{I}_{#1}(#2)}}
\renewcommand{\empty}{\ensuremath{\mbox{empty}}}
\newcommand{\elem}{\ensuremath{\mbox{elem}}}
\renewcommand{\cl}[1]{\ensuremath{\mathscr{#1}}}

\newcommand{\treemodel}[2]{\ensuremath{\mathrm{Treemod}_{#1}(#2)}}

\newcommand{\tm}[2]{\ensuremath{\mathrm{TM}_{#1}(#2)}}
\newcommand{\tmfin}[2]{\ensuremath{\mathrm{TM}^{\text{f}}_{#1}(#2)}}
\newcommand{\tmf}[2]{\ensuremath{\tmfin{r}{d}}}
\newcommand{\mymod}[1]{\ensuremath{\text{Mod}(#1)}}

\newcommand{\tmtree}[2]{\ensuremath{\mathrm{Tree}_{#1}(#2)}}

\newcommand{\mystar}[1]{\ensuremath{{#1}^\star}}
\newcommand{\mytilde}[1]{\ensuremath{\widetilde{#1}}}
\newcommand{\msf}[1]{\ensuremath{\mathsf{#1}}}

\titleformat{\section}[block]{\Large\sc\filcenter}{\thesection.}{5pt}{}
\titleformat{\subsection}[block]{\sc\filcenter}{\thesubsection.}{5pt}{}

\makeatletter
\providecommand{\institute}[1]{
\apptocmd{\@author}{\end{tabular}  \par\smallskip \begin{tabular}[t]{c} #1}{}{} }
\providecommand{\email}[1]{
  \apptocmd{\@author}{\end{tabular}
    \par\smallskip
    \begin{tabular}[t]{c}
    \texttt{#1}}{}{}
}
\makeatother

\title{\tbf{Pseudo-finiteness of arbitrary graphs\\ of bounded shrub-depth}\footnote{This research was supported by the Leverhulme Trust through a Research Project Grant on ``Logical Fractals''.}}
\author{Abhisekh Sankaran}
\institute{\small{Department of Computer Science and Technology,}\\\small{University of Cambridge, UK}}
\email{\small{abhisekh.sankaran@cl.cam.ac.uk}}
\date{}

\begin{document}
\maketitle

\begin{abstract}
We consider classes of arbitrary (finite or infinite) graphs of bounded shrub-depth, specifically the classes $\tm{r}{d}$ of arbitrary graphs that have tree models of height $d$ and $r$ labels. We show that the graphs of $\tm{r}{d}$ are $\mso$-pseudo-finite relative to the class $\tmfin{r}{d}$ of finite graphs of $\tm{r}{d}$; that is, that every $\mso$ sentence true in a graph of $\tm{r}{d}$ is also true in a graph of $\tmfin{r}{d}$. We also show that $\tm{r}{d}$ is closed under ultraproducts and ultraroots. These results have two consequences. The first is that the index of the $\mso[m]$-equivalence relation on graphs of $\tm{r}{d}$ is bounded by a $(d+1)$-fold exponential in $m$. The second is that $\tm{r}{d}$ is exactly the class of all graphs that are $\mso$-pseudo-finite relative to $\tmf{r}{d}$. 
\end{abstract}

\section{Introduction}\label{section:intro}

Pseudo-finite model theory is a branch of model theory that studies the class of finite structures by studying the expansion of this class with infinite structures that have ``finitary'' behaviour. In that, every first order ($\fo$) sentence true in such an infinite structure is also true in a finite structure. These infinite structures are called \emph{pseudo-finite}. The notion  usually is made to also include finite structures which are trivially pseudo-finite. While the study of pseudo-finite structures has been pursued in mathematics since at least the '80s~\cite{field-arithm}, the study of such structures with reference to computer science is more recent, being largely initiated in~\cite{vaananen-pseudo-fin}. The mentioned paper develops the subject as an alternative approach to studying finite models 
as contrasted with finite model theory that typically either studies extensions of FO over all finite structures~\cite{KV92,rosen-weinstein, DHK95}, or studies FO (and its extensions) over restricted classes of finite structures~\cite{ADK06, dawar-pres-under-ext, otto11, HHS15}. The class of all pseudo-finite structures forms an elementary class (that is, definable by an $\fo$ theory) -- indeed it is the class of arbitrary models of the $\fo$ theory of the class of all finite structures. This is an equivalent definition of pseudo-finiteness. Given that most results of $\fo$ model theory relativize to elementary classes, pseudo-finite structures are naturally model-theoretically well-behaved. 

Shrub-depth~\cite{shrub-depth-definitive} is a graph parameter that has been introduced in the context of obtaining algorithmic meta theorems for model checking properties expressible in a well-studied extension of $\fo$ called Monadic Second Order logic ($\mso$), with improved  dependence on the sizes of the input $\mso$ sentences considered as a parameter. In contrast to the usual non-elementariness for the dependence on this parameter for even $\fo$ model checking over the class of all trees~\cite{frick-grohe},  graph classes of bounded shrub-depth admit fixed parameter tractability of $\mso$ model checking with a parameter dependence that is a fixed tower of exponentials, of height proportional to the shrub-depth. These classes are defined using tree models of height $d$ that use at most $r$ labels for some naturals $d$ and $r$, where informally, such a tree model $\tree{t}$ of a graph $G$ is a rooted tree whose leaves are the vertices of $G$, that are assigned labels from the set $\{1, \ldots, r \}$. The presence or absence of an edge between two vertices of $G$ is determined by the labels of these vertices in $\tree{t}$ and the distance between them in $\tree{t}$.  Since its inception, shrub-depth has seen plenty of active research for not just its algorithmic properties, but also its structural and logical aspects~\cite{shrub-depth-definitive,shrub-depth-oum,CF20, GK20}.

In this paper, we study a relativization of the notion of pseudo-finiteness to classes of graphs of bounded shrub-depth, and go further to consider the version of it for $\mso$. That is, instead of considering $\fo$ and all finite graphs in the definition of pseudo-finiteness, we instead consider $\mso$ and the class $\tmf{r}{d}$ of finite graphs that have tree models of height $d$ and $r$ labels. This gives us the class of arbitrary graphs that are \emph{$\mso$-pseudo-finite} \emph{relative to} $\tmf{r}{d}$. The graphs in this class are those for which every $\mso$ sentence true in the graph is true also in some graph of $\tmfin{r}{d}$. Equivalently, these are the  graphs satisfying the $\mso$ theory of $\tmf{r}{d}$. Of central interest to us in the paper is understanding the graphs structurally, akin to this understanding for $\tmf{r}{d}$. Towards this study, we consider \emph{arbitrary} graphs that have tree-models of height $d$ and $r$ labels, where the tree-models could now be infinite. We denote this class of graphs $\tm{r}{d}$. Clearly $\tmf{r}{d}$ is the class of finite graphs of $\tm{r}{d}$. As the first central result of this paper, we show the following.

\begin{theorem}\label{dispthm:tmrd-pseudofin}
For every $m \ge 1$ and every graph $G$ of $\tm{r}{d}$, there exists a graph $H$ in $\tmfin{r}{d}$ such that: (i) $H$ is an induced subgraph of $G$; (ii) $G$ and $H$ agree on all sentences of MSO of quantifier rank at most $m$; and (iii) the size of $H$ is at most a $d$-fold exponential in $m$ and $r$. Thus, in particular, the graphs of $\tm{r}{d}$ are $\mso$-pseudo-finite relative to $\tmfin{r}{d}$.
\end{theorem}

Theorem~\ref{dispthm:tmrd-pseudofin} in fact shows a stronger property than $\mso$-pseudo-finiteness for any infinite graph $G$ of $\tm{r}{d}$; namely that, every $\mso$ sentence true in $G$ is also true in a finite \emph{induced subgraph} of $G$. Theorem~\ref{dispthm:tmrd-pseudofin}  can therefore also be seen as showing a strong form of the classic L\"owenheim-Skolem theorem from model theory for graphs of  $\tm{r}{d}$, by which, not only is it the case that any $\fo$ sentence true in an infinite graph $G$ is true also in a countable induced subgraph of $G$, but also that same holds with $\mso$ instead of $\fo$ and `finite' instead of countable. Further, as a consequence of the bounds provided, we obtain also that the index of the $\mso[m]$-equivalence relation over $\tm{r}{d}$, which relates two graphs of $\tm{r}{d}$ if they agree on all MSO sentences of rank at most $m$, is an elementary function of $m$, indeed a  $(d+1)$-fold exponential in $m$. This is in contrast to the usual non-elementary lower bound for this index even for FO over the class of all finite trees.

Theorem~\ref{dispthm:tmrd-pseudofin} tells us that the graphs of $\tm{r}{d}$ are models of the $\mso$ theory of $\tmfin{r}{d}$. We go further to investigate what other graphs are models of this theory. Before considering infinite models, one first observes that in contrast to the case of the class of all graphs that are ($\fo$-)pseudo-finite (relative to all finite graphs), where the class naturally includes all finite graphs, it is a non-trivial question whether there are any finite graphs not in $\tmf{r}{d}$, that are $\mso$-pseudo-finite (or even $\fo$-pseudo-finite) relative to $\tmf{r}{d}$. It turns out that there are no such graphs. This is because $\tmfin{r}{d}$ has a characterization in terms of a finite number of excluded induced subgraphs~\cite{shrub-depth-definitive}, and therefore $\tmfin{r}{d}$ is axiomatized in the finite by a universal FO sentence  that describes this characterization. This sentence hence belongs to the $\fo$ and $\mso$ theories of $\tmfin{r}{d}$, and indeed axiomatizes both of these in the finite.

Moving onward to the infinite, we now ask what  infinite models other than those in $\tm{r}{d}$ does the $\mso$ theory of $\tmfin{r}{d}$ have. 
This isn't an easy question as such since there aren't many tools to deal with infinite structures for their \emph{$\mso$} properties, as there are for their $\fo$ properties. Since the graphs of $\tm{r}{d}$ are also $\fo$-pseudo-finite, we examine the infinite models of the \emph{FO theory} of $\tmf{r}{d}$. Here we prove the second central result of this paper.
\begin{theorem}\label{dispthm:ultraprodroot-closure}
The class $\tm{r}{d}$ is closed under ultraproducts and ultraroots.
\end{theorem}
The above theorem in conjunction with the fact that $\tm{r}{d}$ is closed under isomorphisms, readily gives us that $\tm{r}{d}$ is an elementary class using a well-known characterization of elementariness under the mentioned closure properties (Theorem~\ref{thm:keisler-shelah-2}). And this inference in conjunction with Theorem~\ref{dispthm:tmrd-pseudofin} gives us the following characterization of $\mso$-pseudo-finiteness (and also $\fo$-pseudo-finiteness) relative to $\tmf{r}{d}$.

\begin{theorem}\label{dispthm:char-pseudofin}
The class $\tm{r}{d}$ is exactly the class of arbitrary models of the $\fo$ theory, and hence also the $\mso$ theory, of $\tmfin{r}{d}$. As a consequence, $\tm{r}{d}$ is characterized over all graphs by the same finite set of excluded finite induced subgraphs that characterizes  $\tmfin{r}{d}$ over finite graphs.
\end{theorem}

The main tool we use for showing Theorem~\ref{dispthm:tmrd-pseudofin} is a version of the Feferman-Vaught  composition  theorem proved in~\cite{FO-MSO-coincide}. 
This version allows  evaluating any $\mso$ sentence $\Phi$ over the disjoint union of an arbitrary family $\mc{F}$ of structures, by examining the truth of an $\fo$ sentence $\alpha_\Phi$ over an \emph{$\mso[m]$-type indicator} $\ind{m}{\mc{F}}$. This  is a  structure over a monadic vocabulary, that contains all the information about the equivalence classes of the $\mso[m]$ relation to which the structures in $\mc{F}$ belong (indeed these classes constitute the mentioned vocabulary). 
The fact that the vocabulary of the type indicator is \emph{monadic} allows us to shrink the structure to under a certain threshold size (that depends on the rank) without any change in the $\fo[q]$ theory in going to the shrunk structure, where $q$ is the considered rank. Using this simple observation, we first prove our results for trees of bounded height noting that trees are after all constructed inductively from forests of lesser height, and the latter lend themselves to using the mentioned composition theorem. Subsequently, the results for trees are transferred to $\tm{r}{d}$ using the $\fo$ interpretability of the latter in the former. The proof above turns out to give elementary sized small models for $\mso$ over $\tmf{r}{d}$ and is a considerably simpler proof than the one for the mentioned result shown in~\cite{shrub-depth-definitive}. For Theorems~\ref{dispthm:ultraprodroot-closure} and~\ref{dispthm:char-pseudofin}, we use a combination of combinatorial and infinitary (compactness based) reasoning to show the results.

The organization of the paper is as follows. In Section~\ref{section:prelims}, we provide the background and notation for the paper; in Section~\ref{section:tmrd-pseudo-fin}, we show Theorem~\ref{dispthm:tmrd-pseudofin} and in  Section~\ref{section:char-pseudo-fin}, we show Theorems~\ref{dispthm:ultraprodroot-closure} and~\ref{dispthm:char-pseudofin}. We present our conclusions in Section~\ref{section:conclusion}.

\vspace{-1pt}

\section{Background}\label{section:prelims}

We assume the reader is familiar with the  terminology in connection with the syntax and semantics $\fo$ and $\mso$. \eat{A vocabulary, typically denoted $\tau$ or $\sigma$, is a finite set of relation symbols. We denote the class of all $\mso$ formulae over $\tau$ as $\mso(\tau)$ and its subclass of $\fo$ formulae as $\fo(\tau)$.} A sequence $x_1, \ldots, x_n$ of $\fo$ variables is denoted $\bar{x}$. An $\mso$ formula $\varphi$ whose free variables are in $\bar{x}$ is denoted $\varphi(\bar{x})$. A sentence is a formula without free variables. The \emph{quantifier rank}, or simply rank, of an $\mso$ formula $\varphi$, denoted $\rank{\varphi}$, is the maximum number of quantifiers (both first order and second order) appearing in any root to leaf path in the parse tree of the formula. We denote by $\mso[m]$ the class of all $\mso$ formulae of rank at most $m$.

A simple, undirected graph is an $\{E\}$-structure in which the binary relation $E$ is interpreted as an irreflexive and symmetric relation. All graphs in the paper are simple and undirected. A \emph{tree} is a connected graph that does not contain any cycles. A $p$-labeled rooted tree $\tree{t}$ is a $\{E, \rootref, P_1, \ldots, P_p\}$-structure such that: (i) the $\{E\}$-reduct of $\tree{t}$ is a tree; (ii) the unary relation symbol $\rootref$ is interpreted as a set consisting of a single element called the \emph{root} of $\tree{t}$, and denoted $\root{\tree{t}}$; and (iii) the $P_i$s for $i \in [p] = \{1, \ldots, p\}$ form a partition of the nodes of $\tree{t}$ (some of the parts could be empty). We shall often call $p$-labeled rooted trees, as simply trees when $p$ is clear from context. We use the standard notions of parent, child, ancestor and descendent in the context of trees. For a tree $\tree{t}$ and a node $v$ of it, the \emph{subtree of $\tree{t}$ rooted at $v$}, denoted $\tree{t}_v$, is the substructure of $\tree{t}$ induced by the descendants of $v$ in $\tree{t}$ (these include $v$), except for the interpretation of the $\rootref$ predicate, which is $\{v\}$ (instead of $\emptyset$). A $p$-labeled rooted forest $\forest{f}$ is a disjoint union of $p$-labeled rooted trees  $\tree{t}_i$  for $i$ belonging to an index set $I$; we then write $\forest{f} = \bigcupdot_{i \in I} \tree{t}_i$ where $\bigcupdot$ denotes disjoint union. A ($p$-labeled rooted) tree $\tree{t}_2$ is said to be a \emph{leaf-hereditary} subtree of a tree $\tree{t}_1$ if the roots of $\tree{t}_2$ and $\tree{t}_1$ are the same, and there is a subset $X$ of nodes of $\tree{t}_1$ deleting the subtrees rooted at which the resulting tree is $\tree{t}_2$. Equivalently, $\tree{t}_2$ is a leaf-hereditary subtree of $\tree{t}_1$ if the roots of $\tree{t}_2$ and $\tree{t}_1$ are the same, $\tree{t}_2$ is a substructure of $\tree{t}_1$, and every leaf of $\tree{t}_2$ is also a leaf of $\tree{t}_1$. 
If $\tree{t}_2$ is a leaf-hereditary subtree of $\tree{t}_1$, then it follows that: (i) for any two leaf nodes of $\tree{t}_2$, the distance between them in $\tree{t}_2$ is the same as the distance between them in $\tree{t}_1$; (ii) if $\forest{f}_i$ is the forest of $p$-labeled rooted trees obtained by removing the root of $\tree{t}_i$ for $i \in \{1, 2\}$, then for every tree $\tree{s}_2$ of $\forest{f}_2$, there exists a tree $\tree{s}_1$ of $\forest{f}_1$ such that $\tree{s}_2$ is a leaf-hereditary subtree of $\tree{s}_1$.
The \emph{height} of a tree $\tree{t}$ is the maximum root to leaf distance in $\tree{t}$.
We denote by $\cl{T}_{d, p}$ the class of  arbitrary (finite or infinite) $p$-labeled rooted trees of height at most $d$. 

\paragraph{Shrub-depth:}\label{section:shrub-depth}
We recall the notion of tree models from~\cite{shrub-depth-definitive} and state it in its extended version for arbitrary cardinality graphs. For $r, d \in \mathbb{N}$ where $\mathbb{N}$ denotes the set of naturals including 0, a \emph{tree model of $r$ labels and height $d$} for a graph $G$ is a pair $(\tree{t}, S)$ where $\tree{t}$ is an $(r+1)$-labeled arbitrary rooted tree of height $d$ and $S \subseteq [r]^2 \times [d]$ is a set called the \emph{signature} of the tree model such that:
\vspace{2pt}\begin{enumerate}
    \item The length of every root to leaf path in $\tree{t}$ is exactly $d$.
    \item The set $V(G)$ is exactly the set of leaves of $\tree{t}$.
    \item Each leaf has a unique label from $[r]$ and all internal nodes are labeled $r+1$.
    \item For any $i, j \in [r]$ and $l \in [d]$, it holds that $(i, j, l) \in S$ if and only if $(j, i, l) \in S$. 
    \item For vertices $u, v \in V(G)$, if $i$ and $j$ are the labels of $u$ and $v$ seen as leaves of $\tree{t}$, and the distance between $u$ and $v$ in $\tree{t}$ is $2l$, then $\{u, v\} \in E(G)$ iff $(i, j, l) \in S$. Observe that the distance between $u$ and $v$ in $\tree{t}$ is an even number as all root to leaf paths are of length $d$, and $l$ is thus the distance between $u$ (or $v$) and the least common ancestor of $u$ and $v$.
\end{enumerate}

\noindent The class of arbitrary tree models of $r$ labels and height $d$ is denoted $\treemodel{r}{d}$ and the class of trees contained in these tree models is denoted $\tmtree{r}{d}$; so $\tmtree{r}{d} = \{\tree{t} \mid (\tree{t}, S) \in \treemodel{r}{d}~\mbox{for some signature}~S \subseteq [r]^2 \times [d]\}$. The class of arbitrary graphs that have tree models in $\treemodel{r}{d}$ is denoted $\tm{r}{d}$, and the class of finite graphs of $\tm{r}{d}$ is denoted $\tmf{r}{d}$. Recalling from~\cite{shrub-depth-definitive}, a class of finite graphs has \emph{shrub-depth} at most $d$ it is a subclass of $\tmf{r}{d}$; then analogously we can define a class of arbitrary graphs to have shrub-depth at most $d$ if it is a subclass of $\tm{r}{d}$. It is easy to see that $\tm{r}{d}$ and $\tmf{r}{d}$ are both hereditary classes, that is, closed under induced subgraphs.

We make some observations from the definition of $\tm{r}{d}$ that we will need in Section~\ref{section:tmrd-pseudo-fin}. We see that for every signature $S \subseteq [r]^2 \times [d]$, there exists a pair $\Xi_S = (\xi_{V, S}(x), \xi_{E, S}(x, y))$ of $\fo$ formulae that when evaluated on a tree $\tree{s}$ in $\tmtree{r}{d}$ produces the graph $G \in \tm{r}{d}$ of which $(\tree{s}, S) \in \treemodel{r}{d}$ is a tree model. Specifically: (i) the formula $\xi_{V, S}(x)$ says that $x$ is a leaf node; and (ii) the formula $\xi_{E, S}(x, y)$ says that $x$ and $y$ are leaves, and for some $(i, j, l) \in S$, it is the case that $P_i(x)$ and $P_j(y)$ are true, and that the distance between $x$ and $y$ is exactly $2l$. The pair $\Xi_{S}$ is called an \emph{$\fo$ interpretation in $\tmtree{r}{d}$}. Thus $\Xi_{S}$ defines a function from  $\tmtree{r}{d}$ to $\tm{r}{d}$, which also we denote as $\Xi_{S}$; so for $G$ and $\tree{s}$ as above $G = \Xi_{S}(\tree{s})$. Then $\tm{r}{d} = \{\Xi_{S}(\tree{s}) \mid (\tree{s}, S) \in \treemodel{r}{d}\}$. 

We now make some important observations about $\Xi_{S}$. Firstly, the rank $q$ of $\Xi_{S}$, defined as the maximum rank of the formulae appearing in it, is such that $q = O(d)$. Next, if $G = \Xi_{S}(\tree{s})$, then for an $\mso$ formula $\varphi$ of rank $m$ in the vocabulary of $G$, there is an $\mso$ formula that we denote $\Xi_{S}(\varphi)$, of rank $m + q$ in the vocabulary of $\tree{s}$ such that the following holds. (This is a special case of a more general result called the \emph{fundamental theorem of interpretations}.)
$$G \models \varphi ~~~\mbox{iff}~~~ \tree{s} \models \Xi_{S}(\varphi)$$. 
As a consequence, if $\tree{s}_1$ and $\tree{s}_2$ are trees of $\tmtree{r}{d}$, then 
$$\tree{s}_1 \mequiv{m+q} \tree{s}_2 ~~~\longrightarrow~~~ \Xi_{S}(\tree{s}_1) \mequiv{m} \Xi_{S}(\tree{s}_2)$$
Finally for $(\tree{s}, S) \in  \treemodel{r}{d}$, if $\tree{s}'$ is a leaf-hereditary subtree of $\tree{s}$, then $(\tree{s}', S) \in  \treemodel{r}{d}$ and $\Xi_{S}(\tree{s}')$ is an induced subgraph of $\Xi_{S}(\tree{s})$.

\paragraph{Feferman-Vaught composition:}\label{section:FV-decomp} 
Let $\mc{L}$ be one of the logics $\fo$ or $\mso$. Given $m \in \mathbb{N}$ and structures $\str{A}$ and $\str{B}$ from a class $\mc{C}$ of structures (say $p$-labeled rooted trees or unlabeled graphs) over a vocabulary $\tau$, we say $\str{A}$ and $\str{B}$ are $\mc{L}[m]$-equivalent, denoted $\str{A} \lequiv{m} \str{B}$, if $\str{A}$ and $\str{B}$ agree on all $\mc{L}$ sentences of rank at most $m$. It is known that the $\lequiv{m}$ relation has finite index~\cite{chang-keisler}. For $\mc{L} = \mso$, we let $\iota(m, \cl{C})$ denote this index (of the $\mequiv{m}$ relation) restricted to $\cl{C}$. 
Let $\mc{F} = (\str{A}_i)_{i \in I}$ be a family of structures of $\cl{C}$ with disjoint universes, indexed by an index set $I$ of an arbitrary cardinality. Let $m \in \mathbb{N}$ and $\tau_{m,  \cl{C}}$ be the relational vocabulary consisting of a distinct unary predicate symbol for each equivalence class of the $\mequiv{m}$ relation over $\cl{C}$, \eat{$T \in \mdelta{m}{\cl{C}}$,} and containing no other predicate symbols. \eat{Note that we are using the same symbol $T$ to denote both an equivalence class of $\mequiv{m}$ over $\cl{C}$, as well as the unary predicate symbol corresponding to the class for the ease of understanding; whether we mean $T$ as a class or a predicate symbol will be clear from context.}  The \emph{$\mso[m]$-type indicator} for the family $\mc{F}$ is now defined as a $\tau_{m,  \cl{C}}$-structure $\ind{m}{\mc{F}}$ such that: (i) the universe of $\ind{m}{\mc{F}}$ is $I$; and (ii) for $T \in \tau_{m, \cl{C}}$ that corresponds to an equivalence class $\delta$ of $\mequiv{m}$ over $\cl{C}$, the interpretation of $T$ in $\ind{m}{\mc{F}}$ is the set $\{i \in I \mid \delta~\mbox{is the}~\mequiv{m}~\mbox{class of}~\str{A}_i~\mbox{in}~\cl{C}\}$. Observe that for each $i \in I$, there is exactly one predicate $T \in \tau_{m,  \cl{C}}$ such that  $i$ is in the interpretation of $T$ in $\ind{m}{\mc{F}}$. We now have the following theorem from~\cite{FO-MSO-coincide}. (This is the  special case of $w = 0$ and $\mathrm{L} = \mso$ in~\cite[Theorem 14]{FO-MSO-coincide}.)

\begin{theorem}[Theorem 14,~\cite{FO-MSO-coincide}]\label{thm:FV-decomposition}
Let $\cl{C}$ be a class of structures over a vocabulary $\tau$. For every $\mso$ sentence $\Phi$ over $\tau$ of rank $m$, there exists an $\fo$ sentence $\alpha_{\Phi}$ over $\tau_{m, \cl{C}}$ such that if $\mc{F} = (\str{A}_i)_{i \in I}$ is a family of structures of $\cl{C}$ with disjoint 

\noindent universes for an index set $I$ of an arbitrary cardinality, then the following holds:
\[
\ind{m}{\mc{F}} \models \alpha_\Phi ~~~\mbox{if, and only if,}~~~ \bigcupdot_{i \in I} \str{A}_i \models \Phi
\]
Further, if $\,\cl{\widetilde{C}}$ is the class of structures of $\cl{C}$ expanded with (all possible interpretations of) $m$ new unary predicate symbols, then the rank of $\alpha_\Phi$ is $O((\iota(m, \cl{\widetilde{C}}))^{m+1})$.
\end{theorem}

\begin{remark}
In~\cite{FO-MSO-coincide}, the result is actually stated for $\Phi$ which does not contain any $\fo$ variables and whose atomic formulae, instead of being the usual atomic formulae (of the form $x_1 = x_2$, $Y(x)$ for an $\mso$ variable $Y$, and $R(x_1, \ldots, x_r)$ where $R \in \tau$ and $x, x_1, \ldots, x_r$ are $\fo$ variables), are instead of the ``second order'' forms $\empty(X)$ and $\elem(X_1, \ldots, X_r, Z)$ where $X, X_1, \ldots, X_r$ are $\mso$ variables and $Z$ is either an $\mso$ variable or a predicate of $\tau$, with $r$ being the arity of $Z$. The semantics for these atomic forms are as suggested by their names \eat{is intuitive: $\empty(X)$ holds for set $P$ if $P$ is empty, and $\elem(X_1, \ldots, X_r, Z)$ holds for $X_i = P_i$ and $Z = Q$ if $|P_i| = 1$ for $1 \leq i \leq r$ and $P_i \times \cdots \times P_r \subseteq Q$.} Every ``usual'' $\mso$ formula can be converted into an equivalent $\mso$ formula over the mentioned second order atomic formulae, without any change of quantifier rank (see~\cite[page 4]{FO-MSO-coincide}). We have hence recalled~\cite[Theorem 14]{FO-MSO-coincide} in the form stated above in Theorem~\ref{thm:FV-decomposition} which features $\Phi$ as a usual $\mso$ formula. \eat{since these are the kinds of formulae which we work with in this paper.}
\end{remark}

We provide here the justification for the last statement of Theorem~\ref{thm:FV-decomposition} which does not appear explicitly in~\cite{FO-MSO-coincide} but is indeed a consequence of the proof of~\cite[Theorem 14]{FO-MSO-coincide}. We refer the reader to~\cite[Section 3.1, pp. 6 -- 9]{FO-MSO-coincide} to find the formulae and other constructions we refer to in our description here. 

We first observe in the proof of Lemma 8 of~\cite{FO-MSO-coincide}, that the ``capping'' constant $C$ for the formula $\alpha$ is simply the rank of $\alpha$, since $\alpha$ is an FO sentence over a monadic vocabulary $\tau_\Sigma$ (and we  also see a similar such result in Lemma~\ref{lemma:monadic-structures}). Then the number $n$  mentioned in the proof is at most $\mbox{rank}(\alpha) \cdot |\tau_\Sigma|$, whereby the rank of $\beta_z$, and hence the rank of $\beta_\alpha$, is at most $\mbox{rank}(\alpha) \cdot |\mbox{vocab}(\alpha)| + 1$, where $\mbox{vocab}(\alpha)$  denotes the vocabulary of $\alpha$, namely $\tau_\Sigma$. Call this observation (*).

We now come to the proof of Theorem 14 of~\cite{FO-MSO-coincide}, and make the following observations about the rank of $\alpha_\Phi$ following the inductive construction of $\alpha_\Phi$ as given in the proof. For the base cases, since $\gamma^{\mathrm{L}}_\Psi(i)$ is a quantifier-free formula for any $\Psi$, we get that if $\Phi := \empty(X)$, then rank of $\alpha_\Phi$ is 1, and if $\Phi := \elem(X_1, \ldots, X_r, Z)$ or $\Phi := \elem(X_1, \ldots, X_r, R)$, then the rank of $\alpha_\Phi$ is 2 since the width $w$ is 0 by our assumption. If $\Phi$ is a Boolean combination of a set of formulae, then $\mbox{rank}(\alpha_\Phi)$ is the maximum of the ranks of the formulae in the mentioned set. We now come to the non-trivial case when $\Phi := \exists X (\Phi')$. 

We see from~\cite[page 9, para 2]{FO-MSO-coincide} that the rank of $\alpha_\Phi$ is the maximum of the ranks of $\alpha_C$ where $\alpha_C$ is obtained from $\beta_{\alpha_{\Phi'}}$ (denoted as simply $\beta$ in the proof as a shorthand) by substituting the atoms $T(i)$ with the quantifier-free formula $\gamma^{\mathrm{L}}_\Psi(i)$ for a suitable $\Psi$. Then $\mbox{rank}(\alpha_\Phi) = \mbox{rank}(\beta_{\alpha_{\Phi'}})$. It follows from (*) above that $\mbox{rank}(\beta_{\alpha_{\Phi'}}) \leq \mbox{rank}(\alpha_{\Phi'}) \cdot |\mbox{vocab}(\alpha_{\Phi'})| + 1$; call this inequality (**).  
Let $\equiv_{r, \mathrm{L}}$ denote the equivalence relation that relates two structures over the same vocabulary iff they agree on all $\mathrm{L}$ sentences (over the vocabulary of the structures) of rank at most $r$. Then the vocabulary of $\alpha_{\Phi'}$ is the set of all equivalence classes of the $\equiv_{q - 1, \mathrm{L}}$ relation over all structures over the vocabulary of the family $F$ (where $F$ is as in the statement of~\cite[Theorem 14]{FO-MSO-coincide}), expanded with (all possible interpretations of) $d$ set predicates, where $d$ is the number of free variables of $\Phi'$ and $q - 1$ is the rank of $\Phi'$. Here we now importantly observe that if the structures of the family $F$ come from a class $\cl{C}$, then it is sufficient to consider just those equivalence classes of the $\equiv_{q - 1, \mathrm{L}}$ relation that are non-empty when restricted to the structures of $\cl{C}$ expanded with $d$ set predicates. Then if $\Phi'$ is a subformula of a rank $m$ $\mathrm{L}$ sentence $\Phi$ over a vocabulary $\tau$ and we are interested only in a given class $\cl{C}$ of $\tau$-structures and expansions of these with set predicates, then the size of the vocabulary of $\alpha_{\Phi'}$ is at most the index of the $\equiv_{m, \mathrm{L}}$ equivalence relation over the class $\widetilde{\cl{C}}$ of structures of $\cl{C}$ expanded with (all possible interpretations of exactly) $m$ set predicates. Applying this observation to (**) iteratively, we then get that if $\mathrm{L} = \mso$ and $\lambda = \iota(m, \cl{\widetilde{C}})$ then

\[
\begin{array}{lll}
\mbox{rank}(\alpha_\Phi) &  \leq & 1 + \lambda \cdot \Big(1 + \lambda \cdot \big( 1 + \ldots \cdot (1 + 2 \cdot \lambda)\big)\Big)\\
     & \leq &  1 + \lambda + \ldots + \lambda^{m-1} + 2 \cdot \lambda^m\\
     & \leq & 2 \cdot \lambda^{m+1}\\
     & = & O((\iota(m, \cl{\widetilde{C}}))^{m+1})
\end{array}
\]
showing the last statement of Theorem~\ref{thm:FV-decomposition}.

\paragraph{Ultraproducts:} Given a family $(\str{A}_i)_{i \in I}$ of structures over a relational vocabulary $\tau$ and an ultrafilter $U$ on the index set $I$, let $\mytilde{\str{A}} = \prod_{i \in I} \str{A}_i$ be the direct (Cartesian) product of the $\str{A}_i$'s. Let $\sim$ be the equivalence relation on the universe of $\mytilde{\str{A}}$ defined as: if $\bar{a} = (a_i)_{i \in I}$ and $\bar{b} = (b_i)_{i \in I}$ are tuples of elements from $\mytilde{\str{A}}$, then $\bar{a} \sim \bar{b}$ if, and only if, $\{i \in I \mid a_i = b_i\} \in U$. Let $[\bar{a}]$ denote the equivalence class of $\bar{a}$ under $\sim$. Then the \emph{ultraproduct} of $(\str{A}_i)_{i \in I}$ is the structure $\mystar{\str{A}} = \prod_{i \in I} \str{A}_i / U$ defined as: (i) the universe of $\mystar{\str{A}}$ is the set $\{ [\bar{a}] \mid \bar{a} \in \mytilde{\str{A}}\}$; (ii) for a $k$-ary relation $R \in \tau$ and tuples $\bar{a}_j = (a_j^i)_{i \in I}$ of $\mytilde{\str{A}}$ for $j \in [k]$, it holds that $\mystar{\str{A}} \models R([\bar{a}_1], \ldots, [\bar{a}_k])$ if, and only if, $\{i \in I \mid \str{A}_i \models R(a_1^i, \ldots, a_k^i)\} \in U$. Given the definition of $\sim$, it can be seen that the presented interpretation of $R$ in $\mystar{\str{A}}$ is well-defined. If $\str{A}_i = \str{A}$ for all $i \in I$, then $\mystar{\str{A}}$ is called the \emph{ultrapower} of $\str{A}$ with respect to $U$, and $\str{A}$ is called the \emph{ultraroot} of $\mystar{\str{A}}$ with respect to $U$. Two well-known theorems concerning the ultraproduct are as below.

\begin{theorem}[{\L}o\'s theorem; Theorem 4.1.9~\cite{chang-keisler}]\label{thm:los-ultraprod}
Let $\mytilde{\str{A}}$ and $\mystar{\str{A}}$ be the direct product and ultraproduct respectively of a family $(\str{A}_i)_{i \in I}$ of structures over a vocabulary $\tau$, with respect to an ultrafilter $U$ over an index set $I$. Then for an $\fo$ formula $\varphi(x_1, \ldots, x_l)$ over $\tau$ and elements $\bar{a}_j = (a_j^i)_{i \in I}$ of 
$\mytilde{\str{A}}$ for $j \in [l]$, 
\[
\mystar{\str{A}} \models \varphi([\bar{a}_1], \ldots, [\bar{a}_l])~~~\mbox{if, and only if,}~~~\{i \in I \mid \str{A}_i \models \varphi(a_1^i, \ldots, a_l^i)\} \in U
\]
\end{theorem}

\begin{theorem}[Theorems 4.1.12 and 6.1.15~\cite{chang-keisler}]\label{thm:keisler-shelah-2}
A class of structures is elementary iff it is closed under isomorphisms, ultraproducts and ultraroots.
\end{theorem}

For a class $\str{C}$ of finite structures and $\mc{L} \in \{\fo, \mso\}$, we denote by $\lth{\mc{C}}$ the class of all $\mc{L}$ sentences that are true in all structures of $\mc{C}$. We say an arbitrary structure $\str{A}$ is \emph{$\mc{L}$-pseudo-finite relative to $\mc{C}$} if $\str{A} \models \lth{\mc{C}}$; that is (equivalently) if every sentence of $\lth{\{\str{A}\}}$ is true in some structure of $\mc{C}$.

\section{Pseudo-finiteness of $\tm{r}{d}$ relative to $\tmf{r}{d}$}\label{section:tmrd-pseudo-fin}

Following are the central results of this section. Define the function $\mathsf{tower}: \mathbb{N} \times \mathbb{N} \rightarrow \mathbb{N}$ as: $\tower{0, n} = n$ and $\tower{d, n} = 2^{\tower{d-1, n}}$. The core technical result is Theorem~\ref{thm:tree-pseudo-fin-index-bound} that shows a relativized $\mso$-pseudo-finiteness theorem for $\cl{T}_{d, p}$. This is then transferred to $\tm{r}{d}$  in Theorem~\ref{thm:pseudo-fin-and-index-bound} using $\fo$ interpretations.

\begin{theorem}\label{thm:pseudo-fin-and-index-bound}
Let $d, r \in \mathbb{N}$ be given. There exists an increasing function $h: \mathbb{N} \rightarrow \mathbb{N}$ such that the following are true for each $m \in \mathbb{N}_+$.  
\vspace{1pt}\begin{enumerate}
    \item  
    For every graph  $G \in \tm{r}{d}$, there exists $H \in \tm{r}{d}$ such that: (i) $H \subseteq G$; (ii) $|H|$ is at most $\tower{d, h(d) \cdot m \cdot (m + \log r)}$, and (ii) $H \mequiv{m} G$. \label{thm:pseudo-fin-and-index-bound:1} 
    \item 
    The index $\iota(m, \tm{r}{d})$ of the $\mequiv{m}$ relation over $\tm{r}{d}$ is at most $\tower{d+1, h(d) \cdot m^2 \cdot (\log r)^2}$.\label{thm:pseudo-fin-and-index-bound:2}
\end{enumerate}
\end{theorem}

\begin{theorem}\label{thm:tree-pseudo-fin-index-bound}
Let $d, p \in \mathbb{N}$ be given. There exists an increasing function $g: \mathbb{N} \rightarrow \mathbb{N}$ such that if $\zeta_{d, p}: \mathbb{N} \times \mathbb{N} \rightarrow \mathbb{N}$ is the function given by $\zeta_{d, p}(n_1, n_2) = \tower{n_2, g(d) \cdot (n_1 + 1) \cdot (n_1 + \log p)}$, then the following are true for each $m \in \mathbb{N}_+$.  
\vspace{1pt}\begin{enumerate}
    \item  
    For every tree $\tree{t} \in \cl{T}_{d, p}$, there exists a leaf-hereditary subtree $\tree{t}'$ of $\tree{t}$ such that: (i) the heights of $\tree{t}'$ and $\tree{t}$ are the same; (ii) $|\tree{t}'|$ is at most $\zeta_{d, p}(m, d)$; and (iii) $\tree{t}' \mequiv{m} \tree{t}$.\label{thm:tree-pseudo-fin-index-bound:1} 
    \item 
    The index $\iota(m, \cl{T}_{d, p})$ of the $\mequiv{m}$ relation over $\cl{T}_{d, p}$ is at most $\zeta_{d, p}(m, d+1)$ if $d \ge 1$, and is $p$ if $d = 0$. \label{thm:tree-pseudo-fin-index-bound:2}
\end{enumerate}
\end{theorem}

\begin{proof}[Proof of Theorem~\ref{thm:pseudo-fin-and-index-bound}]
Since $G \in \tm{r, p}{d}$, from Section~\ref{section:shrub-depth} there exists a tree model $(\tree{s}, S) \in \treemodel{r}{d}$ for $G$ and an $\fo$ interpretation $\Xi_{S}$ such that $G = \Xi_{S}(\tree{s})$. Let $q$ be the rank of $\Xi_S$; we know that $q = O(d)$. Since $\tree{s}$ can be seen as a tree in $\cl{T}_{d, r + 1}$, by Theorem~\ref{thm:tree-pseudo-fin-index-bound}(\ref{thm:tree-pseudo-fin-index-bound:1}), there exists a leaf-hereditary subtree $\tree{s}'$ of $\tree{s}$ such that (i) $|\tree{s}'| \leq \zeta_{d, r + 1}(q + m, d)$, and (ii) $\tree{s}' \mequiv{q + m} \tree{s}$. Since every root-to-leaf path of $\tree{s}'$ is also a root-to-leaf path of $\tree{s}$, and every internal node of $\tree{s}'$ is labeled with the label $r+1$ while the leaf nodes are labeled with labels from $[r]$, we get that $(\tree{s}', S) \in \treemodel{r}{d}$ and that $(\tree{s}', S)$ is a tree model for $H = \Xi_{S}(\tree{s}') \in \tm{r}{d}$. From the properties of the interpretation $\Xi_{S}$ as discussed in Section~\ref{section:shrub-depth}, we infer the following: (i) Since $\tree{s}'$ is a leaf-hereditary subtree of $\tree{s}$, we have $H \subseteq G$; (ii) Since $V(H)$ is the set of leaves of $\tree{s}'$, we have $|H| \leq |\tree{s}'| \leq \zeta_{d, r + 1}(q + m, d) = \tower{d, g(d) \cdot (q + m + 1) \cdot (q + m + \log (r + 1)} \leq  \tower{d, g(d) \cdot q^2 \cdot m \cdot (m + \log r)} \leq \tower{d, h(d) \cdot m \cdot (m + \log r)}$ where $h(d) = c_0 \cdot g(d) \cdot d^2$,  $q \leq c_1 \cdot d $ and $c_0 \ge (c_1)^2$; (iii) Since $\tree{s}' \mequiv{q+m} \tree{s}$, we have $H \mequiv{m} G$. 
\newcommand{\tms}[2]{\ensuremath{\mathrm{TM}^S_{#1}(#2)}}

We now look at the index of the $\mequiv{m}$ relation over $\tm{r}{d}$. For a given signature $S \subseteq [r]^2 \times [d]$, denote $\tms{r}{d}$ denote the subclass of $\tm{r}{d}$ of those graphs that have a tree model $(\tree{t}, S) \in \treemodel{r}{d}$. Then $\Xi_{S}$ is a surjective map from $\tmtree{r}{d}$ to $\tms{r}{d}$. For $q$ as above, since any equivalence class of the $\mequiv{q+m}$ relation over $\tmtree{r}{d}$ gets mapped by $\Xi_{S}$ to a subclass of an equivalence class of the $\mequiv{m}$ relation over $\tms{r}{d}$, we get by the surjectivity of $\Xi_S$ that $i(m, \tms{r}{d}) \leq i(q + m, \tmtree{r}{d}) \leq i(c_1 \cdot d + m, \cl{T}_{d, r + 1})$ for $c_1$ as above. Then $\iota(m, \tm{r}{d}) \leq \sum_{S \subseteq [r]^2 \times [d]} i(m, \tms{r}{d}) \leq 2^{r^2 \cdot d} \cdot i(c_1 \cdot d + m, \cl{T}_{d, r + 1}) \leq 2^{r^2 \cdot d} \cdot \tower{d+1, g(d) \cdot (c_1 \cdot d + m +1) \cdot (c_1 \cdot d + m + \log (r + 1))} \leq \tower{d+1, h(d) \cdot (m \cdot (\log r))^2}$.
\end{proof}

Towards the proof of Theorem~\ref{thm:tree-pseudo-fin-index-bound}, we will require the following lemma that can be shown using a simple {\ef} game argument.

\begin{lemma}\label{lemma:monadic-structures}
Let $\sigma$ be a finite vocabulary consisting of only monadic relation symbols. Let $\str{A}$ be an arbitrary structure over $\sigma$ such that every element of $\str{A}$ is in the interpretation of exactly one predicate of $\sigma$. Let $a$ be a given element of $\str{A}$, and let $q \in \mathbf{N}_+$ and $q_1 = 1 + (q - 1) \cdot |\sigma|$. Then there exists a substructure $\str{B}$ of $\str{A}$ such that: (i) $\str{B}$ contains $a$, (ii) $|\str{B}| \leq q_1$, and (iii) $\str{B} \fequiv{q} \str{A}$.
\end{lemma}

\begin{proof}[Proof of Theorem~\ref{thm:tree-pseudo-fin-index-bound}]
We prove the theorem by induction on $d$. The base case of $d = 0$ is trivial to see by taking $g(0) = 1$. Assume as induction hypothesis that the statement is true for $d - 1$ for $d \ge 1$. 

Consider a tree $\tree{t} \in \cl{T}_{d, p}$ of height equal to $d$. Let $\forest{f} = \bigcupdot_{i \in I} \tree{s}_i$ be the forest of $p$-labeled rooted trees obtained by removing the root of $\tree{t}$. Let $\cl{S}$ be the class of rooted forests whose constituent trees belong to $\cl{T}_{d-1, p}$; so $\forest{f} \in \cl{S}$. Consider now the $\mso[m]$-type indicator $\ind{m}{\mc{F}}$ for the family $\mc{F} = (\tree{s}_i)_{i \in I}$ and the $\mso[m]$ sentence $\Phi$ that axiomatizes the $\mequiv{m}$ equivalence class of $\forest{f}$ (this is known to exist~\cite{libkin}).
By Theorem~\ref{thm:FV-decomposition}, there is an FO sentence $\alpha_\Phi$ over the vocabulary $\tau_{m, \cl{T}_{d-1, p}}$ such that 
\[
\ind{m}{\mc{F}} \models \alpha_\Phi~~\mbox{if, and only if,}~~\forest{f} \models \Phi
\]

Now we know from Theorem~\ref{thm:FV-decomposition} that for any $\mso[m]$ sentence $\Psi$ over the vocabulary of $\cl{T}_{d-1, p}$, the sentence $\alpha_\Psi$ given by the theorem has rank that is $O((\iota(m, \widetilde{\cl{T}}_{d-1, p}))^{m+1})$ where $\widetilde{\cl{T}}_{d-1, p}$ is the expansion of $\cl{T}_{d-1, p}$ with $m$ new unary predicates. Now there is a natural 1-1 correspondence between $\widetilde{\cl{T}}_{d-1, p}$ and $\cl{T}_{d-1, p \cdot 2^m}$, and two $\widetilde{\cl{T}}_{d-1, p}$ structures are $\mso[m]$-equivalent iff their corresponding $\cl{T}_{d-1, p \cdot 2^m}$ structures are. Then the rank of $\alpha_\Psi$ is $O((\iota(m, \cl{T}_{d-1, p \cdot 2^m}))^{m+1})$. Let $c_0$ be a constant such that $\mbox{rank}(\alpha_\Psi) \leq c_0 \cdot (\iota(m, \cl{T}_{d-1, p \cdot 2^m}))^{m+1}$ for all $m \ge m_0$. In fact, as shown in Section~\ref{section:FV-decomp}, one can take $c_0 = 2$ and $m_0 = 1$.

From the analysis above then, we have that $\mbox{rank}(\alpha_\Phi) \leq c_0 \cdot (\iota(m, \cl{T}_{d-1, p \cdot 2^m}))^{m+1}$. Now by induction hypothesis, we see that: (i) if $d - 1 = 0$, then $\iota(m, \cl{T}_{d-1, p \cdot 2^m}) =  p \cdot 2^m$; and (ii) if $d > 1$, then  $\iota(m, \cl{T}_{d-1, p \cdot 2^m}) \leq \zeta_{d-1, p \cdot 2^m}(m, d) = \tower{d, g(d-1) \cdot (m + 1) \cdot (2m + \log p))} \leq \tower{d, 2 \cdot g(d-1) \cdot (m + 1) \cdot (m + \log p)}$. Then the rank of $\alpha_\Phi$ is at most 
\vspace{1pt}\begin{enumerate}
    \item $c_0 \cdot (p \cdot 2^m)^{m+1} \leq  2^{c_0 \cdot (m + 1) \cdot (m + \log p))}$  if $d = 1$, and 
    \item $c_0 \cdot (\tower{d,  2 \cdot g(d-1) \cdot (m + 1) \cdot (m + \log p)})^{m+1} \leq \tower{d, c_0 \cdot 4 \cdot g(d-1) \cdot (m + 1) \cdot (m + \log p)}$, if $d > 1$.
\end{enumerate} 

If $\rho_{d, p}: \mathbb{N} \rightarrow \mathbb{N}$ is the function given by $\rho_{d, p}(m) = \tower{d+1, 4 \cdot c_0 \cdot g(d) \cdot (m + 1) \cdot (m + \log p)}$, then we see that in either case above, the rank of $\alpha_\Phi$ is at most $\rho_{d-1, p}(m)$. We also see that $\tau_{m, \cl{T}_{d-1, p}}$ is the vocabulary which contains one unary predicate symbol for every equivalence class of the $\mequiv{m}$ relation over $\cl{T}_{d-1, p}$, and only those predicates; then $|\tau_{m, \cl{T}_{d-1, p}}| = \iota(m, \cl{T}_{d-1, p})$ which is equal to $p$ if $d - 1 = 0$, and at most $\tower{d, g(d-1) \cdot (m + 1) \cdot (m + \log p))}$ if $d > 1$. Then $\mbox{rank}(\alpha_\Phi) \cdot |\tau_{m, \cl{T}_{d-1, p}}| \leq \rho_{d-1, p}(m) \cdot |\tau_{m, \cl{T}_{d-1, p}}| \leq \tower{d, 5 \cdot c_0 \cdot g(d-1) \cdot (m + 1) \cdot (m + \log p)}$ for all $d \ge 1$.

We observe now that $\ind{m}{\mc{F}}$ is a structure over a finite monadic vocabulary such that each element of its universe is in the interpretation of exactly one predicate in the vocabulary. Let $\tree{s}_{i^*} \in \mc{F}$ for $i^* \in I$ be such that the height of $\tree{s}_{i^*}$ is equal to $d-1$ (there must be such a tree in $\mc{F}$ since height of $\tree{t}$ is equal to $d$). Recall that $I$ is the universe of $\ind{m}{\mc{F}}$. Then by Lemma~\ref{lemma:monadic-structures}, taking $a = i^*$ and $q = \rho_{d-1, p}(m)$ in the lemma, we get that there exists a substructure $\str{B}$ of $\ind{m}{\mc{F}}$ such that (i) $\str{B}$ contains $i^*$, (ii) $|\str{B}| \leq 1 + (q - 1) \cdot |\tau_{m, \cl{T}_{d-1, p}}| \leq q \cdot |\tau_{m, \cl{T}_{d-1, p}}| \leq \tower{d, 5 \cdot c_0 \cdot g(d-1) \cdot (m + 1) \cdot (m + \log p)}$, and (iii) $\str{B} \fequiv{q} \ind{m}{\mc{F}}$. Then $\str{B}$ can be seen as the $\mso[m]$-type indicator $\ind{m}{\mc{F}'}$ of the family $\mc{F}' = (\tree{s}_j)_{j \in I'}$. for a subset $I' \subseteq I$, that contains $i^*$. Then by Theorem~\ref{thm:FV-decomposition}, we have
\[
\ind{m}{\mc{F}'} \models \alpha_\Phi~~\mbox{if, and only if,}~~\bigcupdot_{j \in I'} \tree{s}_j \models \Phi
\]
Since (i) $\ind{m}{\mc{F}'} = \str{B} \fequiv{q} \ind{m}{\mc{F}}$, (ii) $\ind{m}{\mc{F}} \models \alpha_\Phi$, and (iii) $\mbox{rank}(\alpha_\Phi) \leq q$, we have $\ind{m}{\mc{F}'} \models \alpha_\Phi$ and therefore $\bigcupdot_{j \in I'} \tree{s}_j \models \Phi$. Then $\bigcupdot_{j \in I'} \tree{s}_j \mequiv{m} \forest{f}$. 

We now utilize the induction hypothesis for Part (\ref{thm:tree-pseudo-fin-index-bound:1}). Since $\tree{s}_j \in \cl{T}_{d-1, p}$ for all $j \in I'$, there exists a leaf-hereditary subtree $\tree{s}'_j$ of $\tree{s}_j$ such that: (i) the heights of $\tree{s}'_j$ and $\tree{s}_j$ are the same; (ii) $|\tree{s}'_j| \leq \zeta_{d-1, p}(m, d-1)$; and (ii) $\tree{s}'_j \mequiv{m} \tree{s}_j$. Then consider the forest $\forest{f}' = \bigcupdot_{j \in I'} \tree{s}'_j$. Since disjoint union satisfies the Feferman-Vaught composition property, we get that $\forest{f}' \mequiv{m} \bigcupdot_{j \in I'} \tree{s}_j \mequiv{m} \forest{f}$. Further, we have $1 + |\forest{f}'| = 1 + \sum_{j \in I'} |\tree{s}'_j| \leq  1 + |I'| \cdot \max_{j \in I'} |\tree{s}'_j| \leq 1 +  \tower{d, 5 \cdot c_0 \cdot g(d-1) \cdot (m + 1) \cdot (m + \log p)} \cdot \tower{d-1, g(d-1) \cdot (m + 1) \cdot (m + \log p)} \leq \tower{d, 6 \cdot c_0 \cdot g(d-1) \cdot (m + 1) \cdot (m + \log p)}$ for all $d \ge 1$. Let $\tree{t}'$ be the leaf-hereditary subtree of $\tree{t}$ obtained by removing all child subtrees (that are the unrooted versions of) $\tree{s}_j$ for $j \notin I'$, and replacing $\tree{s}_j$ with $\tree{s}'_j$ (again the replacement being for the unrooted versions of the trees) for each $j \in I'$ in $\tree{t}$. Then $\forest{f'}$ is indeed the forest of rooted trees obtained by deleting the root of $\tree{t}'$. Observe that $\forest{f}'$ contains the tree $\tree{s}'_{i^*}$ (since $i^* \in I'$) whose height is the same as that of $\tree{s}_{i^*}$ whose height in turn is equal to $d-1$; then $\tree{t}'$ has height equal to $d$ which is the height of $\tree{t}$. Further, since $\forest{f}' \mequiv{m} \forest{f}$, it is easy to verify that $\tree{t}' \mequiv{m} \tree{t}$. Finally, $|\tree{t}'| = |\forest{f}'|+ 1 \leq \tower{d, 6 \cdot c_0 \cdot g(d-1) \cdot (m + 1) \cdot (m + \log p)}$ for all $d \ge 1$. 

We now observe that the existence of $\tree{t}'$ as above for every tree $\tree{t}$ in $\cl{T}_{d, p}$ implies that $\iota(m, \cl{T}_{d, p})$ is at most the number of structures of $\cl{T}_{d, p}$ whose size is at most $\tower{d, 6 \cdot c_0 \cdot g(d-1) \cdot (m + 1) \cdot (m + \log p)}$. Since the number of structures of $\cl{T}_{d, p}$ with universe size at most $\mu$ for any number $\mu$ is at most $\mu \cdot 2^{\mu \cdot (\mu + \log p)} \leq 2^{3 \mu^2}$ if $\log p \leq \mu$, we get, by taking $\mu = \tower{d, 6 \cdot c_0 \cdot g(d-1) \cdot (m + 1) \cdot (m + \log p)}$, that $\iota(m, \cl{T}_{d, p}) \leq \tower{d+1, 14 \cdot c_0 \cdot g(d-1) \cdot (m + 1) \cdot (m + \log p)}$. Then defining $g(d) = 14 \cdot c_0 \cdot g(d-1)$, we see that both parts of the present theorem are true for $d$. This completes the induction and hence the proof.
\end{proof}

\begin{remark}\label{remark:alpha-phi-properties}
By the same reasoning as in the proof of Theorem~\ref{thm:tree-pseudo-fin-index-bound}, it follows that if $\cl{C}$ in Theorem~\ref{thm:FV-decomposition} is taken to be $\cl{T}_{d, p}$, then for the $\mso[m]$ sentence $\Phi$ as considered in Theorem~\ref{thm:FV-decomposition} (that is, an arbitrary $\mso$ sentence of rank $m$ over the vocabulary of $\cl{T}_{d, p}$), the sentence $\alpha_\Phi$ is such that: (i) the vocabulary $\tau_{m, \cl{T}_{d, p}}$ of $\alpha_\Phi$ has size at most $\zeta_{d, p}(m, d+1)$; (ii) the rank of $\alpha_\Phi$ is at most $\rho_{d, p}(m)$; and (iii) $\mbox{rank}(\alpha_\Phi) \cdot |\tau_{m, \cl{T}_{d, p}}| \leq \zeta_{d+1, p}(m, d+1)$.
\end{remark}

A corollary of Theorem~\ref{dispthm:tmrd-pseudofin} that will be useful for us in the next section is the following.  
\begin{corollary}\label{cor:equal-theories}
For every $r, d \ge 0$, it holds that $\foth{\tm{r}{d}} = \foth{\tmf{r}{d}}$ and $\msoth{\tm{r}{d}} = \msoth{\tmf{r}{d}}$.
\end{corollary}
\begin{proof}
We show the result for $\mso$; the proof for $\fo$ is analogous. That $\msoth{\tm{r}{d}} \subseteq \msoth{\tmf{r}{d}}$ is obvious; for the other direction, suppose it is false. Then for some $G \in \tm{r}{d}$ and $\varphi \in \msoth{\tmf{r}{d}}$, it holds that $G \models \neg \varphi$. Then $H \models \neg \varphi$ for some $H \in \tmf{r}{d}$ by Theorem~\ref{thm:pseudo-fin-and-index-bound} which is a contradiction.
\end{proof}

\section{Characterizing pseudo-finiteness relative to $\tmf{r}{d}$}\label{section:char-pseudo-fin}

In this section, we complement the results of Section~\ref{section:tmrd-pseudo-fin} by showing that the graphs of $\tm{r}{d}$ are exactly those that are $\mso$-pseudo-finite relative to $\tmf{r}{d}$. At the heart of this result is the following theorem. 
\begin{theorem}\label{thm:ultraprodroot-closure}
The class $\tm{r}{d}$ is closed under ultraproducts and ultraroots.
\end{theorem}
\begin{proof}
\tbf{Ultraproducts:} Let $(G_i)_{i \in I}$ be an indexed family of graphs in $\tm{r}{d}$ and $U$ be an ultrafilter on $I$. Let $\mystar{G} = \prod_{i \in I} G_i / U$ be the ultraproduct of the $G_i$'s with respect to $U$. For $k \in I$, let $(\tree{t}_k, S_k) \in \treemodel{r}{d}$ be a tree model for $G_k$. Define the pair $(\mystar{\tree{t}}, \mystar{S})$ such that:
\[
 \mystar{\tree{t}} = \prod_{i \in I} \tree{t}_i / U~~~~; ~~~~~ (i, j, l) \in \mystar{S} ~~~\mbox{iff}~~~ \{k \in I \mid (i, j, l) \in S_k\} \in U
\]
We show below that $(\mystar{\tree{t}}, \mystar{S}) \in \treemodel{r}{d}$, and then show that $(\mystar{\tree{t}}, \mystar{S})$ is a tree model for an isomorphic copy of $\mystar{G}$ to complete the proof of this part of the theorem. 

Consider the FO sentences below; recall that the vocabulary of the trees $\tree{t}_k$ is $\tau = \{E, \mathsf{root}\} \cup \{P_i \mid i \in [r+1]\}$.

\begin{tabular}{llp{.8\textwidth}}
    $\chi$ & $:=$ & $\bigwedge_{i \in [6]} \chi_i$ \\
    $\chi_1$ & $:=$ & ``There is exactly one vertex satisfying $\root{\cdot}$ (that will be called `root' below)''\\
    $\chi_2$ & $:=$ & ``Every vertex has a path of length at most~$d$~to the root, and there is exactly one such path''\\
    $\chi_3$ & $:=$ & ``If~$x$~and~y~are adjacent and~$x$~is at a distance of~$l$~from the root, then~$y$~is at a distance of~$l_1$~from the root, where~$l_1 \in \{l-1, l+1\}$ and $l \in \{0, \ldots, d\}$''\\
    $\chi_4$ & $:=$ & ``There is no vertex  that is at a distance of at most $d-1$ from the root, and that has degree 1''\\
    $\chi_5$ & $:=$ & ``Every vertex $x$ at a distance of at most~$d-1$~from the root satisfies~$P_{r+1}(x)\wedge \bigwedge_{i \in [r]} \neg P_i(x)$''\\
    $\chi_6$ & $:=$ & ``Every vertex $x$ at a distance equal to~$d$~from the root satisfies $(\bigvee_{i \in [r]} P_i(x)) \wedge (\bigwedge_{i, j \in [r+1], i \neq j} \neg (P_i(x) \wedge P_j(x)))$''\\
\end{tabular}\\

We see that $\{k \in I \mid \tree{t}_k \models \chi\} = I \in U$; then $\mystar{\tree{t}} \models \chi$. So $\mystar{\tree{t}}$ is a (possibly infinite) $(r+1)$-labeled rooted tree of height at most $d$, whose leaf nodes are labeled using labels from $[r]$ and whose internal nodes are labeled with the label $r+1$; in other words $\mystar{\tree{t}} \in \tmtree{r}{d}$ and $(\mystar{\tree{t}}, \mystar{S}) \in \treemodel{r}{d}$. 


Let now $\mytilde{\tree{t}} = \prod_{i \in I} \tree{t}_i$ and $\mytilde{G} = \prod_{i \in I} G_i$, and let $\sim_{\tree{t}}$ and $\sim_G$ resp. be the equivalence relations on $\mytilde{\tree{t}}$ and $\mytilde{G}$ such that the equivalence classes of these relations are the universes of $\mystar{\tree{t}}$ and $\mystar{G}$. For a tuple $\bar{a}$ in $\mytilde{\tree{t}}$, resp. $\mytilde{G}$, we denote by $[\bar{a}]_{\tree{t}}$, resp. $[\bar{a}]_G$ the equivalence class of $\bar{a}$ under $\sim_{\tree{t}}$, resp. $\sim_G$. So  $V(\mystar{\tree{t}}) = \{ [\bar{a}]_{\tree{t}} \mid \bar{a} \in \mytilde{\tree{t}} \}$ and  $V(\mystar{G}) = \{ [\bar{a}]_G \mid \bar{a} \in \mytilde{G}\}$. We recall from the definition of the ultraproduct in Section~\ref{section:prelims}, that for tuples $\bar{a} = (a_k)_{k \in I}$ and $\bar{b} = (b_k)_{k \in I}$ where $a_k, b_k \in V(\tree{t}_k)$, we have
that $\bar{a} \sim_{\tree{t}} \bar{b}$ iff $\{ k \in I \mid a_k = b_k \} \in U$. Likewise for $\sim_G$. Call a tuple $\bar{a} = (a_k)_{k \in I}$ of $\mytilde{\tree{t}}$ \emph{genuinely leaf}, resp. \emph{genuinely non-leaf}, if $a_k$ is a leaf node, resp. non-leaf node, of $\tree{t}_k$ for all $k \in I$. We observe the following. 

\begin{lemma}\label{lemma:reps}
For every tuple $\bar{a}$ of $\mytilde{\tree{t}}$, there is a tuple $\bar{b}$ of $\mytilde{\tree{t}}$ such that: (i) $[\bar{a}]_{\tree{t}} = [\bar{b}]_{\tree{t}}$, and (ii) if $[\bar{a}]_{\tree{t}}$ is a leaf node, resp. non-leaf node, of $\mystar{\tree{t}}$, then $\bar{b}$ is a genuinely leaf, resp. genuinely non-leaf, tuple of $\mytilde{\tree{t}}$.  
\end{lemma}
\begin{proof}
Let $\bar{a} = (a_k)_{k \in I}$, and let $\beta(x)$ be the FO formula that asserts that $x$ is at a distance of exactly $d$ from the root.  Then $[\bar{a}]_{\tree{t}}$ is a leaf node of $\mystar{\tree{t}}$ iff $\mystar{\tree{t}} \models \beta([\bar{a}]_{\tree{t}})$ iff $\{ k \in I \mid \tree{t}_k \models \beta(a_k)\} \in U$ iff $X = \{ k \in I \mid a_k~\mbox{is a leaf node of}~\tree{t}_k \} \in U$. 

Let $\bar{b} = (b_k)_{k \in I} \in V(\mytilde{\tree{t}})$ be defined as follows. 
\begin{enumerate}
    \item If $[\bar{a}]_{\tree{t}}$ is a leaf node of $\mystar{\tree{t}}$, then let $c_k$ be a leaf node of $\tree{t}_k$ for $k \notin X$. Define $b_k = a_k$ if $k \in X$, and $b_k = c_k$ otherwise. We observe that $\{ k \in I \mid a_k = b_k \} = X$. Since $X \in U$ (from above), we have $[\bar{a}]_{\tree{t}} = [\bar{b}]_{\tree{t}}$.
    \item If $[\bar{a}]_{\tree{t}}$ is a non-leaf node of $\mystar{\tree{t}}$, then let $c_k$ be a non-leaf node of $\tree{t}_k$ for $k \in X$. Define $b_k = a_k$ if $k \notin X$, and $b_k = c_k$ otherwise. We observe that $Y = \{ k \in I \mid a_k = b_k \} = \bar{X} = I \setminus X$. Since $X \notin U$ (from above), we have $\bar{X} \in U$ since $U$ is an ultrafilter, and therefore $Y \in U$. Then $[\bar{a}]_{\tree{t}} = [\bar{b}]_{\tree{t}}$.
\end{enumerate}
Thus in both cases $[\bar{a}]_{\tree{t}} = [\bar{b}]_{\tree{t}}$. Further, $\bar{b}$ is a genuinely leaf  if $[\bar{a}]_{\tree{t}}$ is a leaf node of $\mystar{\tree{t}}$, otherwise 
$\bar{b}$ is genuinely non-leaf.
\end{proof}

Using the above lemma, we now argue that the leaf nodes of $\mystar{\tree{t}}$ are in 1-1 correspondence with the vertices of $\mystar{G}$. We firstly observe that since the nodes of $G_k$ are exactly the leaf nodes of $\tree{t}_k$ for $k \in I$, the set $V(\mytilde{G})$ is exactly the set of the genuinely leaf tuples of $V(\mytilde{\tree{t}})$. Now for two leaf nodes $[\bar{a}]_{\tree{t}}$ and $[\bar{b}]_{\tree{t}}$ of $\mystar{\tree{t}}$, let $\bar{a}'$ and $\bar{b}'$ be genuinely leaf tuples given by Lemma~\ref{lemma:reps} such that $\bar{a} \sim_{\tree{t}} \bar{a}'$ and $\bar{b} \sim_{\tree{t}} \bar{b}'$. Then  $\bar{a} \sim_{\tree{t}} \bar{b}$ iff $\bar{a}' \sim_{\tree{t}} \bar{b}'$ iff $\{ k \in I \mid a_k' = b_k' \} \in U$ iff $\bar{a}' \sim_G \bar{b}'$. Let $f$ be the function from the leaf nodes of $\mystar{\tree{t}}$ to $V(\mystar{G})$ such that $f([\bar{a}]_{\sim_{\tree{t}}}) = [\bar{a}']_{\sim_{G}}$  From the equivalences just mentioned, we verify that $f$ is indeed a bijection.
Now for $i, j \in [r]$ and $l \in [d]$, let $\chi_{i, j, l}(x, y)$ be the FO formula as below:
\[
\chi_{i, j, l}(x, y) := U_i(x) \wedge U_j(y) \wedge ``\text{There is a path of length}~2l~\text{between}~x~\text{and}~y"
\]
We now show that for two nodes $[\bar{a}]_G$ and $[\bar{b}]_G$ of $\mystar{G}$ such that $\mystar{\tree{t}} \models \chi_{i, j, l}([\bar{a}]_{\tree{t}}, [\bar{b}]_{\tree{t}})$, it holds that $\mystar{G} \models E([\bar{a}]_G, [\bar{b}]_G)$ iff $(i, j, l) \in \mystar{S}$. This would show that  $(\mystar{\tree{t}}, \mystar{S})$ is a tree-model for an isomorphic copy of $G$, and hence that $G \in \tm{r}{d}$ since $\tm{r}{d}$ is closed under isomorphisms. Recall from the outset that $\mystar{S} = \{ (i', j', l') \mid \{ k \in I \mid (i', j', l') \in S_k\} \in U \}$.

Let $\bar{a} = (a_k)_{k \in I}$ and $\bar{b} =  (b_k)_{k \in I}$. Since $\mystar{\tree{t}} \models \chi_{i, j, l}([\bar{a}]_{\tree{t}}, [\bar{b}]_{\tree{t}})$, we have that $Z = \{ k \in I \mid \tree{t}_k \models \chi_{i, j, l}(a_k, b_k)\} \in U$.  We now obtain the following equivalences. We recall that $(\tree{t}_k, S_k)$ is a tree model for $G_k$ for $k \in I$, and that $U$ is a filter and is hence closed under finite intersections and taking supersets of its members.

\[\def\arraystretch{1.3}
\begin{array}{llll}\baselinestretch
     & \mystar{G} \models  E([\bar{a}]_{G}, [\bar{b}]_G) \\
     \leftrightarrow & \{ k \in I \mid G_k \models E(a_k, b_k)\} \in U \\
     \leftrightarrow & \{ k \in I \mid G_k \models E(a_k, b_k)\} \cap Z \in U\\ 
     \leftrightarrow & \{ k \in I \mid G_k \models E(a_k, b_k) ~\mbox{and}~ \tree{t}_k \models \chi_{i, j, l}(a_k, b_k)\} \in U\\
     \leftrightarrow & \{ k \in I \mid \tree{t}_k \models \chi_{i, j, l}(a_k, b_k) ~\mbox{and}~ (i, j, l) \in S_k\} \in U\\
     \leftrightarrow & Z \cap \{k \in I \mid (i, j, l) \in S_k\} \in U\\
     \leftrightarrow & \{k \in I \mid (i, j, l) \in S_k\} \in U\\
     \leftrightarrow & (i, j, l) \in \mystar{S}
\end{array}
\]

\vspace{2pt} \noindent \tbf{Ultraroots:} 
Let $H$  be an arbitrary graph, and $\mystar{H}$ be an ultrapower of $H$ with respect to an ultrafilter $U$ on an index set $I$.  
Suppose that $\mystar{H}$ belongs to $\tm{r}{d}$ and that $(\mystar{\tree{t}}, \mystar{S}) \in \treemodel{r}{d}$ is a tree-model of $\mystar{H}$. Let $\mytilde{H} = \prod_{i \in I} H$ and let $\sim_H$ be the equivalence relation on $V(\mytilde{H})$ defined as: for  tuples $\bar{d}$ and $\bar{e}$ in $V(\mytilde{H})$, as $\bar{d} \sim_H \bar{e}$ iff $\{i \in I \mid d_i = e_i\} \in U$. Let $[\bar{d}]_{H}$ denote the equivalence class of $\bar{d}$ under $\sim_H$. We know that 
$V(\mystar{H}) = \{ [\bar{d}]_H \mid \bar{d}~\text{is a tuple from}~V(\mytilde{H})\}$. 

Let $X \subseteq V(\mystar{H})$ be defined as $X = \{ [\bar{a}^+]_H \mid a \in V(H)\}$ where 
$\bar{a}^+ = (a, a, a, \ldots)$ is a vertex of $V(\mytilde{H})$. We observe that for distinct vertices $a$ and $b$ of $H$, the vertices $[\bar{a}^+]_H$ and $[\bar{b}^+]_H$ must be distinct as well since $\bar{a}^+ \not\sim_H \bar{b}^+$. We further observe that 
$$ 
\mystar{H} \models E([\bar{a}^+]_H, [\bar{b}^+]_H) ~~~\mbox{iff}~~~ \{i \in I \mid H \models E(a_i, b_i)\} \in U ~~~\mbox{iff}~~~ H \models E(a, b)
$$
since $a_i = a$ and $b_i = b$ for all $I$; so that the map $f: V(H) \rightarrow V(\mystar{H})$ defined as $f(a) = [\bar{a}^+]_H$ is an isomorphic embedding of $H$ into $\mystar{H}$.

We recall that the vertices of $\mystar{H}$ are exactly the leaf nodes of $\mystar{\tree{t}}$. Consider the subtree $\tree{t}$ of $\mystar{\tree{t}}$ induced by $X$ and all the ancestors in $\mystar{\tree{t}}$ of the nodes of $X$. The leaves of $\tree{t}$ are in bijection with the nodes of $H$ via the map $f$. We now have the following equivalences. Below $\msf{lab}([\bar{a}^+]_H; \mystar{\tree{t}})$ denotes the label in $[r]$ of $[\bar{a}^+]_H$ in $\mystar{\tree{t}}$, and $\msf{dist}([\bar{a}^+]_H, [\bar{b}^+]_H; \mystar{\tree{t}})$ denotes the distance between $[\bar{a}^+]_H$ and $[\bar{b}^+]_H$ in $\mystar{\tree{t}}$.
\[\def\arraystretch{1.3}
\begin{array}{ll}
     & H \models E(a, b) \\
    \leftrightarrow & \mystar{H} \models E([\bar{a}^+]_H, [\bar{b}^+]_H) \\
    \leftrightarrow & \mbox{for some}~ i, j \in [r]~\mbox{and}~l \in [d], ~\mbox{it holds that}\\
    & \msf{lab}([\bar{a}^+]_H; \mystar{\tree{t}}) = i, ~\msf{lab}([\bar{b}^+]_H; \mystar{\tree{t}}) = j, ~\msf{dist}([\bar{a}^+]_H, [\bar{b}^+]_H; \mystar{\tree{t}}) = 2l ~\mbox{and}~ (i, j, l) \in \mystar{S}\\
    \leftrightarrow &  \mbox{for some}~ i, j \in [r]~\mbox{and}~l \in [d], ~\mbox{it holds that}\\
    & \msf{lab}([\bar{a}^+]_H; \tree{t}) = i, ~\msf{lab}([\bar{b}^+]_H; \tree{t}) = j, ~\msf{dist}([\bar{a}^+]_H, [\bar{b}^+]_H; \tree{t}) = 2l ~\mbox{and}~ (i, j, l) \in \mystar{S}\\
\end{array}
\]
Taking $S = \mystar{S}$, the equivalences show that $(\tree{t}, S) \in \treemodel{r}{d}$ is indeed a tree-model of an isomorphic copy of $H$. Since $\tm{r}{d}$ is closed under isomorphisms, $H$ belongs to $\tm{r}{d}$.
\end{proof}


Towards proving our desired characterization, we would need another closure property of $\tm{r}{d}$ as shown in Proposition~\ref{prop:another-closure}. We say that a graph class $\cl{C}$ is \emph{closed under membership of finite induced subgraphs} if whenever for a graph $G$ it is the case that all of its finite induced subgraphs are in $\cl{C}$, it holds that $G$ is also in $\cl{C}$. We recall the following well-known result from classical model theory before presenting our result.

\begin{lemma}[Proposition 5.2.2, ref.~\cite{chang-keisler}]\label{lemma:exis-amalgam}
Let $\mf{A}$ and $\mf{B}$ be structures such that every existential sentence that is also true in $\mf{B}$ is true in $\mf{A}$. Then $\mf{B}$
is embeddable in an elementary extension of $\mf{A}$.
\end{lemma}

\begin{proposition}\label{prop:another-closure}
If $\mc{C}$ is a hereditary elementary class, then it is closed under membership of finite induced subgraphs.
\end{proposition}
\begin{proof}
Let $T$ be an FO theory that defines $\mc{C}$. Let $G$ be a simple, undirected graph such that every finite induced subgraph of $G$ belongs to $\mc{C}$. Let $Z$ be the set of all existential FO sentences true in $G$. We show that the theory $T \cup Z$ is satisfiable in a graph $H$. Assuming this to be true, we have by Lemma~\ref{lemma:exis-amalgam}, that there is an elementary extension $H'$ of $H$ within which $G$ embeds isomorphically via an embedding say $f$. Since $H \models T$, we have $H' \models T$ and therefore $H' \in \mc{C}$ since $T$ defines $\mc{C}$. Since $\mc{C}$ is a hereditary class, we have that the image of $G$ under $f$, and therefore $G$, belongs to $\mc{C}$. 

We now show that $T \cup Z$ is satisfiable. For if not, then by Compactness theorem, we have $T \cup Z_1$ is unsatisfiable for a finite subset $Z_1$ of $Z$, and since $Z$ is closed under conjunctions,  we have that $T \cup \{\varphi\}$ is unsatisfiable where $\varphi \in Z$ is the conjunction of the sentences of $Z_1$. Let $Z_1 = \{\psi_i \mid 1 \leq i \leq n\}$ for some $n \ge 1$. Let $A_i \subseteq V(G)$ be the set of witnesses in $G$, to the existential quantifiers in $\psi_i$ (we can always assume $\psi_i$ to be in prenex normal form). Let $A = \bigcup_{i \in [n]} A_i$ and let $G' = G[A] \subseteq G$ be the subgraph of $G$ induced by $A$. We now observe the following:
\begin{enumerate}
    \item $G'$ is a finite induced subgraph of $G$; then $G' \in \mc{C}$ and hence $G'  \models T$.
    \item $G' \models \psi_i$ for $i \in [n]$ because: (i) $G'$ is an extension of $G[A_i]$ (the subgraph of $G$ induced by $A_i$); (ii) the graph $G[A_i]$ models $\psi_i$ and; (iii) $\psi_i$ is preserved under extensions owing to being an existential sentence. Then $G' \models \varphi$.
\end{enumerate}
Then $G' \models T \cup \{\varphi\}$. This is a contradiction with our earlier inference that $T \cup \{\varphi\}$ is unsatisfiable.
\end{proof}

With the above result, we can now obtain the characterization promised at the outset of this section. For a ($\fo/\mso$) theory $T$, let $\mymod{T}$ denote the set of arbitrary models of $T$.

\begin{theorem}\label{thm:char-pseudo-fin}
The following are true for $\mc{L} \in \{\fo, \mso\}$.
\begin{enumerate}
    \item $\tm{r}{d}$ is exactly the class of arbitrary models of $\lth{\tmf{r}{d}}$.\label{thm:char-1}
    \item $\tm{r}{d}$ is characterized over arbitrary graphs by the same finite set of excluded finite induced subgraphs known from~\cite[Theorem 3.10]{shrub-depth-definitive}, that characterizes $\tmfin{r}{d}$ over all finite graphs.\label{thm:char-2}
    \item There exists a universal $\fo$ sentence that axiomatizes $\lth{\tmf{r}{d}}$ over arbitrary graphs.\label{thm:char-3}
\end{enumerate}
\end{theorem}
\begin{proof}
(\ref{thm:char-1}): Now since $\tm{r}{d}$ is isomorphism-closed, we have by Theorems~\ref{thm:ultraprodroot-closure} and~\ref{thm:keisler-shelah-2} that $\tm{r}{d}$ is an elementary class. If $T$ is a theory that defines $\tm{r}{d}$ ($\supseteq \tmf{r}{d}$), then $T \subseteq \foth{\tmf{r}{d}} = \foth{\tm{r}{d}} \subseteq  \msoth{\tm{r}{d}} = \msoth{\tmf{r}{d}}$ by Corollary~\ref{cor:equal-theories}. Then $\tm{r}{d} \subseteq$ $\mymod{\msoth{\tmf{r}{d}}}  \subseteq $  $ \mymod{\foth{\tmf{r}{d}}} \subseteq \mymod{T} = \tm{r}{d}$.

(\ref{thm:char-2}): Let $\mc{F}$ be the finite set of excluded finite induced subgraphs  from~\cite{shrub-depth-definitive} that defines $\tmfin{r}{d}$ in the finite. Let $\varphi$ be the FO sentence that is the conjunction of the negations of the existential closures of the atomic diagrams of the graphs of $\mc{F}$. Then an arbitrary graph $G$ satisfies $\varphi$ iff $G$ excludes the graphs of $\mc{F}$ as induced subgraphs. Then $\varphi$ defines $\tmfin{r}{d}$ over all finite graphs and hence belongs to $\foth{\tmfin{r}{d}}$. Since $\tm{r}{d}$ is defined by $\foth{\tmf{r}{d}}$, every graph of $\tm{r}{d}$ satisfies $\varphi$, and hence excludes the graphs of $\mc{F}$ as induced subgraphs.

Suppose $G$ excludes the graphs of $\mc{F}$ as induced subgraphs; then $G \models \varphi$. Observe that $\varphi$ is equivalent to a universal sentence and is hence hereditary, so that in particular, all finite induced subgraphs of $G$ model $\varphi$. Then all finite induced subgraphs of $G$ belong to $\tmf{r}{d}$ since $\varphi$ defines $\tmfin{r}{d}$ in the finite. Since $\tm{r}{d}$ is a hereditary elementary class, by Proposition~\ref{prop:another-closure}, we get that $G$ belongs to $\tm{r}{d}$.

(\ref{thm:char-3}): The sentence $\varphi$ above, that is equivalent to a universal FO sentence, defines $\tm{r}{d}$ and hence axiomatizes $\lth{\tmf{r}{d}}$ for $\mc{L} \in \{\fo, \mso\}$ over arbitrary graphs.
\end{proof}

\section{Conclusion}\label{section:conclusion}
In this paper, we studied the notion of $\mso$-pseudo-finiteness relative to the class $\tmf{r}{d}$ of finite graphs that have finite tree models of height $d$ and $r$ labels. As our main result, we showed that the class of \emph{arbitrary} graphs that have (arbitrary) tree models of height $d$ and $r$ labels is exactly the class of all graphs that are $\mso$-pseudo-finite relative to $\tmf{r}{d}$. As consequences, we obtained that $\tm{r}{d}$ is characterized over all graphs by the same finite set of excluded finite induced subgraphs known from~\cite{shrub-depth-definitive} to characterize $\tmf{r}{d}$ in the finite, and that the index of the $\mequiv{m}$ relation over $\tm{r}{d}$ is bounded by a $(d+1)$-fold exponential function in $m$.


An interesting consequence of our results is that they allow for transferring results back and forth between $\tm{r}{d}$ and $\tmf{r}{d}$. For instance, one can lift the equivalence of  $\mso$ and $\fo$ over $\tmf{r}{d}$ shown in~\cite{shrub-depth-FO-equals-MSO}, to $\tm{r}{d}$: if an $\mso$ sentence $\varphi$ is equivalent to an $\fo$ sentence $\psi$ over $\tmf{r}{d}$, then `$\varphi \leftrightarrow \psi$' belongs to $\msoth{\tmf{r}{d}}$ and therefore also to $\msoth{\tm{r}{d}}$ by Corollary~\ref{cor:equal-theories}, showing the equivalence of $\varphi$ and $\psi$ over $\tm{r}{d}$. Conversely, classical model-theoretic results like the {\lt} theorem that are true over $\tm{r}{d}$ (since it is an elementary class), can be seen to relativize to $\tmf{r}{d}$. 

For future work, we would like to understand the interaction between the tree model structure of $\tm{r}{d}$ graphs and other structural features of these graphs known from model theory. For instance, whether a graph of $\tm{r}{d}$ is isomorphic to (and not just embeddable in) some  ultraproduct of its finite induced subgraphs. We would also like to explore further the mentioned two-way transfer of results within $\tm{r}{d}$ between the finite and the infinite, including algorithmic meta theorems over graphs of $\tm{r}{d}$ that are finitely presentable.
\vspace{-2pt}

\bibliographystyle{plain}
\bibliography{refs}

\begin{thebibliography}{10}

\bibitem{dawar-pres-under-ext}
Albert Atserias, Anuj Dawar, and Martin Grohe.
\newblock Preservation under extensions on well-behaved finite structures.
\newblock {\em SIAM J. Comput.}, 38(4):1364--1381, 2008.

\bibitem{ADK06}
Albert Atserias, Anuj Dawar, and Phokion~G. Kolaitis.
\newblock On preservation under homomorphisms and unions of conjunctive
  queries.
\newblock {\em J. {ACM}}, 53(2):208--237, 2006.

\bibitem{chang-keisler}
Chen~Chung Chang and H~Jerome Keisler.
\newblock {\em Model theory}.
\newblock Elsevier, 1990.

\bibitem{CF20}
Yijia Chen and J{\"o}rg Flum.
\newblock {FO-Definability of Shrub-Depth}.
\newblock In {\em CSL '20}, LIPIcs, pages 15:1--15:16. Schloss
  Dagstuhl--Leibniz-Zentrum fuer Informatik, 2020.

\bibitem{DHK95}
Anuj Dawar, Lauri Hella, and Phokion~G Kolaitis.
\newblock Implicit definability and infinitary logic in finite model theory.
\newblock In {\em International Colloquium on Automata, Languages, and
  Programming}, pages 624--635. Springer, 1995.

\bibitem{shrub-depth-definitive}
Patrice~Ossona de~Mendez, Jan Obdr{\v{z}}{\'a}lek, Jaroslav
  Ne{\v{s}}et{\v{r}}il, Petr Hlin{\v{e}}n{\`y}, and Robert Ganian.
\newblock Shrub-depth: Capturing height of dense graphs.
\newblock {\em Logical Methods in Computer Science}, 15, 2019.

\bibitem{FO-MSO-coincide}
Michael Elberfeld, Martin Grohe, and Till Tantau.
\newblock Where first-order and monadic second-order logic coincide.
\newblock In {\em {LICS} 2012, Croatia, June 25-28, 2012}, pages 265--274,
  2012.

\bibitem{frick-grohe}
Markus Frick and Martin Grohe.
\newblock The complexity of first-order and monadic second-order logic
  revisited.
\newblock {\em Ann. Pure Appl. Log.}, 130(1-3):3--31, 2004.

\bibitem{field-arithm}
Michael~D Fried and Moshe Jarden.
\newblock {\em Field arithmetic}.
\newblock Springer-Verlag, Heidelberg, 1986.

\bibitem{shrub-depth-FO-equals-MSO}
Jakub Gajarsky and Petr Hlinen{\'{y}}.
\newblock Kernelizing {MSO} properties of trees of fixed height, and some
  consequences.
\newblock {\em Log. Meth. Comp. Sci.}, 11(19):1--26, 2015.

\bibitem{GK20}
Jakub Gajarsk{\`y} and Stephan Kreutzer.
\newblock Computing shrub-depth decompositions.
\newblock In {\em 37th International Symposium on Theoretical Aspects of
  Computer Science (STACS 2020)}. Schloss Dagstuhl-Leibniz-Zentrum f{\"u}r
  Informatik, 2020.

\bibitem{HHS15}
Frederik Harwath, Lucas Heimberg, and Nicole Schweikardt.
\newblock Preservation and decomposition theorems for bounded degree
  structures.
\newblock {\em Log. Methods Comput. Sci.}, 11(4), 2015.

\bibitem{KV92}
Phokion~G. Kolaitis and Moshe~Y. Vardi.
\newblock Fixpoint logic vs. infinitary logic in finite-model theory.
\newblock In {\em LICS '92}, pages 46--57. {IEEE} Computer Society, 1992.

\bibitem{shrub-depth-oum}
O.~Kwon, R.~McCarty, S.~il~Oum, and P.~Wollan.
\newblock Obstructions for bounded shrub-depth and rank-depth.
\newblock {\em J. Comb. Theory, Ser. {B}}, 149:76--91, 2021.

\bibitem{libkin}
Leonid Libkin.
\newblock {\em Elements of finite model theory}.
\newblock Springer Science \& Business Media, 2013.

\bibitem{otto11}
Martin Otto.
\newblock Model theoretic methods for fragments of {FO} and special classes of
  (finite) structures.
\newblock {\em Finite and algorithmic model theory}, 379:271--341, 2011.

\bibitem{rosen-weinstein}
Eric Rosen and Scott Weinstein.
\newblock Preservation theorems in finite model theory.
\newblock In {\em International Workshop on Logic and Computational
  Complexity}, pages 480--502. Springer, 1994.

\bibitem{vaananen-pseudo-fin}
Jouko V{\"a}{\"a}n{\"a}nen.
\newblock Pseudo-finite model theory.
\newblock {\em Mat. Contemp}, 24(8th):169--183, 2003.

\end{thebibliography}

\end{document}